\documentclass{paper}

\usepackage{amssymb}
\usepackage{amscd}
\usepackage{amsfonts}
\usepackage{amssymb}
\usepackage{cite}
\usepackage{amsmath}
\usepackage{epsfig}
\usepackage{graphicx, subfigure}
\usepackage{color, graphics}
\usepackage{pdfsync}
\usepackage{hyperref}
\usepackage[all,cmtip]{xy}
\usepackage{setspace}
 \newtheorem{thm}{Theorem}
 \newtheorem{prop}{Proposition}
 \newtheorem{lem}{Lemma}

 \newtheorem{con}{Conjecture}
 \newtheorem{rem}{Remark}
 \numberwithin{equation}{section}

\renewcommand{\le}{\leqslant}
\renewcommand{\ge}{\geqslant}
\renewcommand{\setminus}{\smallsetminus}

\allowdisplaybreaks

\newcommand{\R}{\mathbb{ R}}

\newcommand{\g}{\mathfrak{ g}}

\DeclareMathOperator{\pr}{pr}

\DeclareMathOperator{\diag}{diag}

\DeclareMathOperator{\ddim}{ddim}
\DeclareMathOperator{\dind}{dind}

\def\diag{\mathrm{diag}}

\def\px1{p_{x_1}}
\def\px2{p_{x_2}}
\def\pu1{p_{u_1}}

\begin{document}

\title{Integrability of homogeneous exact magnetic flows on spheres}

\author{Vladimir Dragovi\'c, Borislav Gaji\'c, and Bo\v zidar Jovanovi\'c}

\maketitle

\noindent{\small Department of Mathematical Sciences, The University of Texas at Dallas,  USA \& Mathematical Institute SANU, Belgrade, Serbia
\footnote{{\sc email:} Vladimir.Dragovic@utdallas.edu}}

\noindent{\small Mathematical Institute SANU, Belgrade,
Serbia\footnote{{\sc email:} gajab@mi.sanu.ac.rs}$^,$\footnote{{\sc email:} bozaj@mi.sanu.ac.rs}}

\begin{abstract} We consider motion of a material point placed in a constant homogeneous magnetic field in $\R^n$ and also motion restricted to the sphere $S^{n-1}$.
While there is an obvious integrability of the  magnetic system in $\R^n$, the integrability  of the system restricted to the sphere $S^{n-1}$ is highly non-trivial. We prove complete integrability of the obtained restricted magnetic systems for $n\le 6$. The first integrals of motion of the  magnetic flows on the spheres $S^{n-1}$, for $n=5$ and $n=6$, are polynomials of the degree $1$, $2$, and $3$ in momenta.
We  prove noncommutative integrability of the obtained magnetic flows for any  $n\ge 7$ when the systems allow a reduction to the cases with $n\le 6$.  We conjecture that the restricted magnetic systems on $S^{n-1}$ are integrable for all $n$.
\footnote{{\sc msc:} 37J35, 53D25. {\sc keywords:} magnetic geodesic flows;  Liouville integrability; noncommutative integrability; Dirac magnetic Poisson bracket; gauge Noether symmetries.}
\end{abstract}

\noindent{\it Dedicated to  Academician V. V. Kozlov on the occasion of his 75-th anniversary.}

\section{Introduction}

We consider motion of a material point of mass $m$ in a constant homogeneous magnetic field in $\R^n$
given by the two-form
\begin{equation*}
\mathbf F=s\sum_{i<j}\kappa_{ij} d\gamma_i\wedge d\gamma_j,
\end{equation*}
and restricted to the sphere
\[
S^{n-1}=\{(\gamma_1, \dots, \gamma_n)\in\R^n\,\vert\, \langle \gamma,\gamma\rangle=\sum_{j=1}^n\gamma_j^2=1\}\subset \R^n.
\]
Here $\kappa =(\kappa_{ij})\in \mathrm{so}(n)$ is a fixed skew-symmetric matrix and $s\in \mathbb R\setminus \{0\}$ is a real parameter representing the charge (with the minus sign) of a material point for $n=3$ (see Remark \ref{3D}).
In \cite{DGJ2023} we obtained these magnetic systems as a reduction of the nonholonomic problem of rolling of a ball with the gyroscope without slipping and twisting over a plane and over a sphere in $\R^n$, where the inertia operator of the system ball + gyroscope is proportional to the identity operator.  We called these nonholonomic systems \emph{the generalized Demchenko cases}. We also proved integrability of the magnetic systems  on $S^2$ and $S^3$,
which correspond to  $n=3$ and $n=4$ respectively, and we performed explicit integrations of the equations of motion of these two systems in elliptic functions  in  \cite{DGJ2023}. See also Theorem  \ref{stara}  below.

As a side benefit, we consider motion of a material point placed in a constant homogeneous magnetic field in $\R^n$.
The system is a direct generalization of the well known 3-dimensional case: for even $n$ generic motions are quasi-periodic over tori of dimension $n/2$ and for odd $n$
generic motions are quasi-periodic over cylinders with the base being tori of dimension $[n/2]$ (see Theorem \ref{ocigledna}).

In this paper we prove complete integrability of the magnetic systems on spheres $S^4$ and $S^5$, corresponding to $n=5$ and $n=6$ respectively, for any $\kappa$  in Theorem \ref{glavna}. The first integrals of motion for these two magnetic systems are polynomials of the  degree $1$, $2$, and $3$ in momenta. The first intergrals that are polynomials linear in momenta are examples of the gauge Noether first integrals (see e.g.  \cite{CS1981}). We also  prove  noncommutative integrability of the obtained magnetic systems for any  $n\ge 7$ when a system allows a reduction to the cases with $n\le 6$ (Theorems \ref{integrabilni}, \ref{integrabilni2}, \ref{integrabilni3}) --- in particular, in the simplest case when
$\mathbf F=s\kappa_{12}d\gamma_1\wedge d\gamma_2$ (item (i), Theorem \ref{integrabilni}) and  for even $n$, when $\mathbf F$ is the magnetic field of the standard contact structure on $S^{n-1}$  (item (i) of Theorem \ref{integrabilni2}).  Thus, it is natural to conjecture that restricted on $S^{n-1}$,  magnetic systems are also integrable for all $n$  and $\kappa$
 (Conjecture \ref{hipoteza}).

Study of dynamics in magnetic fields is one of many areas in which V .V. Kozlov has had a profound impact, see e.g. \cite{Ko1997, KP, VoKo}.

\section{Magnetic flows}

\subsection{Lagrangian systems with magnetic forces}

Let $(Q,\mathbf G)$ be a Riemannian manifold, where $\mathbf G$ is a Riemannian metric. Let
$\mathbf A$ be a one-form on $Q$ and $\mathbf F=d\mathbf A$.
We consider a Lagrangian system $(Q, L_1)$, where the Lagrangian $L_1$ is given by
\[
L_1(q,\dot q)=\frac12 \mathbf G(\dot q,\dot q)+\mathbf A(\dot q)-V(q).
\]

A path $q(t)$ is a \emph{motion of the natural mechanical system} $(Q,L_1)$ if it satisfies the Lagrange-d'Alembert equations
\[
\delta L_1=\big(\frac{\partial L_1}{\partial q}-\frac{d}{dt}\frac{\partial L_1}{\partial \dot q},\delta q\big)=0, \quad \text{for all} \quad \delta q\in T_q Q.
\]
The  Lagrange-d'Alembert equations are equivalent to
\begin{equation}\label{LagMag}
\delta L=\big(\frac{\partial L}{\partial q}-\frac{d}{dt}\frac{\partial L}{\partial \dot q},\delta q\big)=\mathbf F(\dot q,\delta q), \quad \text{for all} \quad \delta q\in T_q Q,
\end{equation}
where $L$ is the part of the Lagrangian $L_1$ without the term linear in velocities
\[
L(q,\dot q)=\frac12\mathbf G(\dot q,\dot q)-V(q).
\]

One can  consider a more general
class of systems \eqref{LagMag}, where an additional force is given as a two-form $\mathbf F$ which is closed, but not need to be exact (see \cite{N1982}). These forces are usually referred as  \emph{magnetic forces}.

\subsection{Hamiltonian description}

Let $(q_1,\dots,q_n)$ be local coordinates on $Q$ in which the Lagrangian and the magnetic two-form take the form:
\begin{align*}
L(q,\dot q)= \frac{1}{2}\sum  g_{ij}(q) \dot q_i\dot q_j-V(q), \qquad \mathbf F=\sum_{i<j} f_{ij}(q) d q_i\wedge d q_j.
\end{align*}

Let $(p_1,\dots,p_n,q_1,\dots,q_n)$ be
the canonical coordinates  of the cotangent bundle $T^*Q$,
and consider the usual Legendre transformation of $L$
\[
p_i=\partial L/\partial \dot q_i=\sum_j g_{ij}\dot q_j, \qquad H=\sum_i p_i\dot q_i-L\vert_{\dot q_i=\sum_j g^{ij}p_j}=\frac12 \sum g^{ij} p_ip_j+V(q),
\]
where $\{g^{ij}\}$ is the inverse of the metric matrix $\{g_{ij}\}$.
Even when the magnetic term is exact $\mathbf F=d\mathbf A$ and we have the Lagrangian $L_1$, it is more convenient to  work with the standard Hamiltonian
\[
H(q,p)=\frac12\mathbf G^{-1}(p,p)+V(q),
\]
rather then the Hamiltonian function that would have linear terms in momenta, obtained as the Legendre transformetion of $L_1$.

In the canonical coordinates, the equations of motion \eqref{LagMag} take the form
\begin{equation}
\label{ham}\dot q_i=\frac{\partial H}{\partial p_i}=\sum_{j=1}^n g^{ij} p_j,\qquad
\dot p_i=-\frac{\partial H}{\partial q_i}+\sum_{j=1}^n f_{ij}(q) \frac{\partial H}{\partial p_j}.
\end{equation}

Let $\rho\colon T^*Q\to Q$ be the canonical projection, $\Omega$ be the canonical symplectic form on $T^*Q$,
$\Omega=dp_1\wedge dq_1+\dots+dp_n\wedge dq_n$, and $z=(q,p)$.
The equations of motion \eqref{ham} can be rewritten in the Hamiltonian form with respect to the \emph{twisted symplectic form} $\Omega+\rho^*\mathbf F$  (see \cite{N1982}):
\begin{equation}\label{red:eq}
\dot z=X_{H}(z), \qquad i_{X_{H}(z)}(\Omega+\rho^*\mathbf F)=-dH(z).
\end{equation}

The problem of integrability
 of magnetic flows has been extesively studied and several integrable cases are known, see e.g.
 \cite{AS2020, BK2017, BJ2008, T2016, MSY2008, S2002}.
We recall a class of integrable magnetic flows on adjoint orbits of compact Lie groups obtained in \cite{BJ2008}, as  an example relevant to this discussion.

\subsection{Magnetic flows on adjoint orbits}

Let $G$ be a compact Lie group,  $\g=T_eG$ be its Lie algebra, and $\langle \cdot,\cdot\rangle$ be an invariant scalar product on $\g$.
Consider the adjoint orbit through an element $a\in\g$:
\[
\mathcal O(a)=\{x=Ad_g(a)=g\cdot a\cdot g^{-1}\,\vert\, g\in G\}.
\]
The adjoint orbit is the homogeneous space
$G/G_a$, where $G_a$ is the isotropy group of $a$.
Since $G$ is a compact connected Lie group, $G_a$ is also connected. Let $\pi\colon G\to G/G_a$ be the natural projection and $\g_a$ be the isotropy algebra of $a$:
$$
\g_a=\{\xi\in\mathfrak g, \, [\xi,a]=0\}=T_e G_a.
$$
The tangent space of $O(a)$ at $\pi(e)$ can be naturally identified by the orthogonal complement to $\g_a$:
\[
T_{\pi(e)}O(a)\cong \g_a^\perp=\{\xi\in\g\, \vert\, \langle \xi,\g_a\rangle=0\}=[a,\mathfrak g].
\]
By definition, the Kirillov-Konstant symplectic form $\Omega_{KK}$ on $O(a)$
 is a $G$--invariant form, given at the point $\pi(e)\in G/G_a$ by
\begin{equation}
\Omega_{KK}(\xi_1,\xi_2)\vert_{\pi(e)}=
-\langle a, [\xi_1,\xi_2]\rangle, \quad   \xi_1,\xi_2\in \g_a^\perp.
\label{Kirillov}
\end{equation}

Consider now a  natural mechanical system on $O(a)$ with the kinetic energy given by the $G$-invariant metric induced by an  invariant scalar product  $\langle \cdot,\cdot\rangle$ and the potential function $V(x)=-\langle b,x\rangle$ under the influence of the magnetic force field given
by $s \Omega_{KK}$, where $s$ is a real non-zero parameter. The Hamiltonian of this system is of the form (see \cite{BJ2008})
\begin{equation}\label{ham-pen}
H(x,p)=\frac12 \langle [x,p],[x,p]\rangle- \langle
b,x\rangle,
\end{equation}
where the cotangent bundle of the orbit is considered as a subspace of the product of Lie algebras:
$$
T^*\mathcal O(a)\subset T^*\mathfrak g \cong \mathfrak g \times \mathfrak g(x,p).
$$
The equations of motion of the magnetic flow, in redundant
variables $(x,p)$ are
\begin{eqnarray}
\dot x=[x,[p,x]], \qquad \dot p=[p,[p,x]] + s [x,p]+b-\pr_{\g_x} b.\label{pp}
\end{eqnarray}

For a generic $b\in\g$, the magnetic system is Liouville integrable.  In the case of the magnetic geodesic flow, when $b=0$, the system is completely integrable
in a non-commutative sense (see \cite{BJ2008}).
As an example, consider the Lie group $SO(3)$. The Lie algebra
$\mathrm{so}(3)$ is isomorphic to the Euclidean space
$\mathbb R^3$ with the bracket operation being the standard  vector product; see the identification \eqref{ident} given below.  From now on, $\langle \cdot,\cdot\rangle$ denotes the standard scalar product in $\R^n$.
The adjoint orbits are spheres. For the unit sphere, the phase space is
\begin{equation*}
T^*S^2=\{(\gamma,p)\in \mathbb R^6 \,\vert\,
\phi_1=\langle \gamma,\gamma\rangle =1, \, \phi_2=\langle \gamma,p\rangle =0\}.
\end{equation*}
The Hamiltonian \eqref{ham-pen} and the system \eqref{pp} reduce to the Hamiltonian of the magnetic spherical pendulum
$$
H=\frac12 \langle p,p\rangle - \langle b, \gamma \rangle,
$$
with the equations of motion
\begin{equation}\label{pendulum}
\dot\gamma=p, \quad \dot p=s \gamma\times p+ b-\big(\langle p,p\rangle+\langle b,\gamma\rangle\big)\gamma.
\end{equation}

 The magnetic term  represents the Lorentz force $s \gamma\times p$ of the magnetic field $\vec B=\gamma$ acting on a material point with the unit mass and the charge $-s$.
The magnetic field $\vec B(\gamma)=\gamma$ in $\R^3$ is known as a magnetic monopole ($\mathrm{div}(\vec B)=3\ne 0$,  compare to the homogeneous field in Remark \ref{3D}).
The system \eqref{pendulum} is completely
integrable and the first integral of motion is
$\langle b, \gamma \times p+s \gamma \rangle$.
In the case when $b=0$, the system \eqref{pendulum}  is non-commutatively completely integrable. There is preservation of the magnetic momentum mapping
$\Phi=\gamma \times p+ s\gamma$ and the trajectories are small circles on the unit sphere.  For a motion
with the unit velocity (magnetic geodesic line), the radius of the circle is equal to
$
r_s={\arctan}({1}/{\vert s \vert}).
$
As $\vert s \vert$ tends to infinity, $r_s$ tends to
zero.  As $s$ tends to zero,  $r_s$ tends to
${\pi}/2$.

\begin{rem}
It is interesting to recall that magnetic monopoles are still not observed in nature. However, we see that
the magnetic term $s\Omega_{KK}$ naturally appears in classical mechanics in the $SO(2)$-reduction of the rigid body motion around a fixed point from $SO(3)$ to $S^2$ for a non-zero value of the integral representing the projection of the angular momentum to the vertical axe, see \cite{NS1981, AKN}.
Therefore, we get integrability of the magnetic pendulums in the cases of the reductions of the Lagrange top and the Kowalevski top.
The above magnetic pendulum \eqref{pendulum} corresponds to the particular case of the Lagrange top when the rigid body is homogeneous.
The reduction of the Euler case represents the magnetic geodesic flow on the sphere $S^2$ endowed with the Poisson sphere metric.
\end{rem}

\section{Motion of a material point in a homogeneous magnetic field}

Consider motion of a material point of mass $m$   in a homogeneous magnetic field in $\R^n$ given by the constant
magnetic two-form $\mathbf F$
\begin{equation}\label{polje}
\mathbf F=s\sum_{i<j}\kappa_{ij} d\gamma_i\wedge d\gamma_j,
\end{equation}
where $\kappa\in \mathrm{so}(n)$ is a fixed skew-symmetric matrix and $s\ne 0$ is a real parameter.
Note that the form $\mathbf F$ is exact:
\begin{equation}\label{AG}
\mathbf F=d\mathbf A^\Gamma, \qquad \mathbf A^\Gamma=\frac{s}{2}\sum_{ij} \kappa_{ij}(\gamma_i+\Gamma_i)d\gamma_j,
\end{equation}
where $\Gamma=(\Gamma_1,\dots,\Gamma_n)\in\R^n$ is an arbitrary fixed vector.

Let $\Omega$ be the canonical symplectic form on $\R^{2n}(\gamma,p)$:
\[
\Omega=dp_1\wedge d\gamma_1+\dots+dp_n\wedge d\gamma_n.
\]
For simplicity, we consider $\mathbf F$ as a two-form on $\R^{2n}$ as well.
The magnetic flow on $(\R^{2n},\Omega+\mathbf F)$ with the Hamiltonian
\begin{equation}\label{HAM}
H(\gamma,p)=\frac{1}{2m}\langle p,p\rangle
\end{equation}
is given by
\begin{align}\label{Rn}
\dot\gamma= \frac{\partial H}{\partial p}= \frac{1}{m} p, \qquad
\dot p   =  -\frac{\partial H}{\partial \gamma}+ s\kappa\big(\frac{\partial H}{\partial p}\big)=\frac{s}{m}\kappa p.
\end{align}

\begin{figure}[h]
{\centering
{\includegraphics[width=7cm]{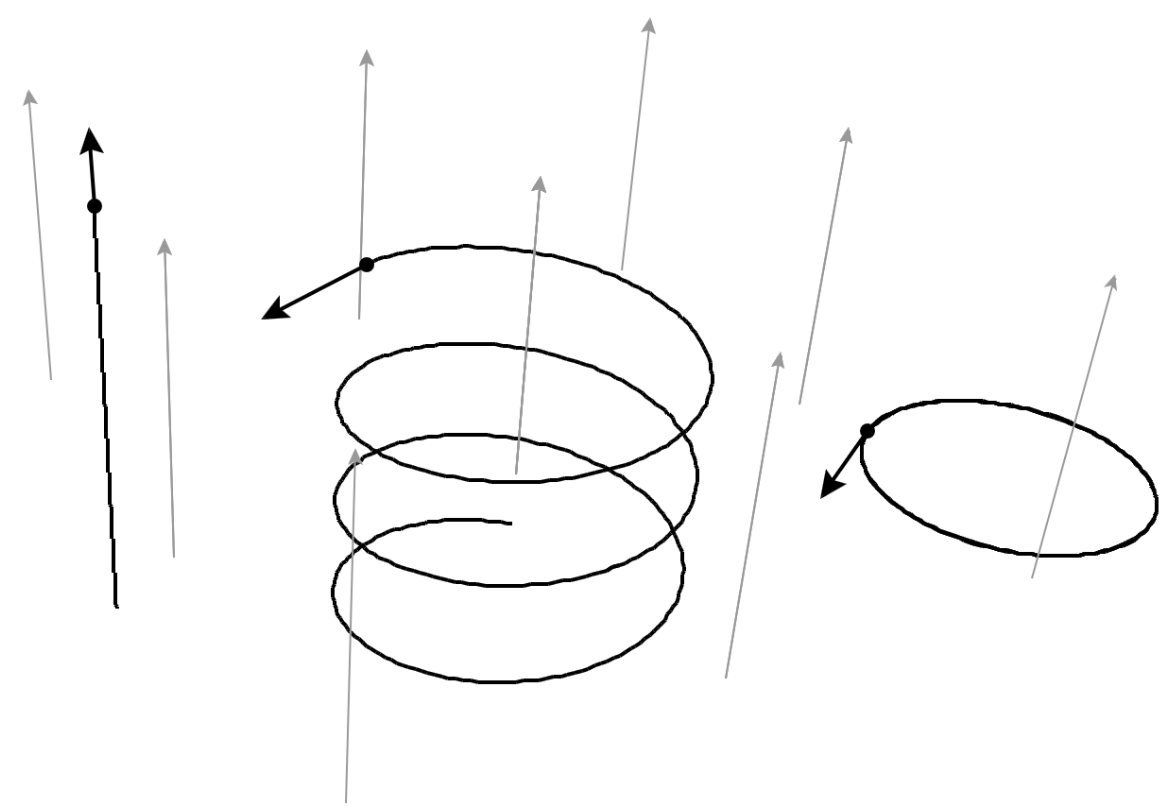}}
\caption{Three types of trajectories in the homogeneous magnetic field in $\R^3$. The gray lines represent the magnetic field. } \label{helikoid}
}
\end{figure}

\begin{rem}\label{3D}
For $n=3$, we aply the standard identification $\R^3\cong \mathrm{so}(3)$,
\begin{equation}\label{ident}
\vec\kappa=(k_1,k_2,k_3) \longleftrightarrow \kappa=\begin{pmatrix}
0 & -k_3 & k_2 \\
k_3 & 0 & -k_1  \\
-k_2 & k_1 & 0
\end{pmatrix}.
\end{equation}
We get a usual form of the Lorentz force
\[
\dot p=\frac{s}{m}\kappa p=\frac{s}{m}\vec\kappa \times p,
\]
 where $\vec\kappa=\vec B$ represents the standard homogeneous magnetic field and $-s$ a charge of a material point of mass $m$.
A trajectory of the system for $\kappa\ne 0$ is  a line directed to $\vec\kappa$,  a circle orthogonal to $\vec\kappa$, or a cylindrical helix
in the direction of $\vec\kappa$ (see Fig. \ref{helikoid}). The projection of a trajectory to the plane orthogonal to $\vec \kappa$ is the circle of the \emph{Larmor radius} $r={\vert p_\perp\vert}/{\vert s\vec\kappa\vert}$ and
the period of rotation   $T={2\pi m}/{\vert s\vec\kappa\vert}$  (see e.g. \cite{F}).
\end{rem}

We have not found the system \eqref{Rn} to be studied before  in the literature.
A lack of interest in this system seems curious, because it is quite natural and  simple.  Without loss of generality,
by taking a suitable orthonormal basis $[\mathbf e_1,\dots,\mathbf e_n]$ of $\R^n$,  we can assume that
the magnetic form \eqref{polje} takes the form
\begin{equation}\label{polje*}
\mathbf F=s(\kappa_{12} d\gamma_1\wedge d\gamma_2+\kappa_{34} d\gamma_3\wedge d\gamma_4 +\dots +\kappa_{2[n/2]-1,2[n/2]}d\gamma_{2[n/2]-1}\wedge d\gamma_{2[n/2]}).
\end{equation}

 If $n$ is even, the system \eqref{Rn} decouples on $k=n/2$ magnetic systems
\begin{align}
\label{res1} &{\dot\gamma}_{2i-1}=\frac{1}{m}p_{2i-1},\qquad\qquad {\dot p}_{2i-1}=\frac{s}{m}\kappa_{2i-1,2i} p_{2i},\\
\label{res2} &{\dot\gamma}_{2i} =\frac{1}{m}p_{2i},\qquad\qquad\qquad  {\dot p}_{2i}=-\frac{s}{m}\kappa_{2i-1,2i} p_{2i-1}.
\end{align}
 They are Hamiltonian equations in $\R^4(\gamma_{2i-1},\gamma_{2i},p_{2i-1},p_{2i})$ with respect to the twisted symplectic forms
\[
\Omega_{2i-1,2i}+\mathbf F_{2i-1,2i}=dp_{2i-1}\wedge g\gamma_{2i-1}+ dp_{2i}\wedge d\gamma_{2i}+s\kappa_{2i-1,2i} d\gamma_{2i-1}\wedge d\gamma_{2i}
\]
and the Hamiltonians
\[
H_{2i-1,2i}=\frac{1}{2m}(p_{2i-1}^2+p_{2i}^2).
\]

If $n=2k+1$ is odd, along with the $k$ systems listed above, there is  an additional  system of one degree of freedom on $\R^2(\gamma_n,p_n)$ with the standard symplectic form and the Hamiltonian
\[
H_n=\frac{1}{2m}p_n^2,
\]
which generates a uniform motion:
\[
\dot\gamma_n=\frac{1}{m}p_n,\qquad\qquad {\dot p}_n=0.
\]

The systems  \eqref{res1}, \eqref{res2} can be easily integrated in terms of linear and trigonometric functions. The projection of the trajectories to $\R^2(\gamma_{2i-1},\gamma_{2i})$-planes are straight lines for $\kappa_{2i-1,2i}= 0$, or circles
\begin{align*}
& \gamma_{2i-1}(t)=r_{2i-1,2i}\cos(\omega_{2i-1,2i}t+\varphi_{2i-1,2i})+c_{2i-1}, \\
&\gamma_{2i}(t)=r_{2i-1,2i}\cos(\omega_{2i-1,2i}t+\varphi_{2i-1,2i})+c_{2i}, \qquad \omega_{2i-1,2i}=-\frac{s\kappa_{2i-1,2i}}{m}
\end{align*}
of the Larmor radii
\[
r_{2i-1,2i}=\frac{\sqrt{p^2_{2i-1}+p^2_{2i}}}{\vert s\kappa_{2i-1,2i}\vert}=\frac{m\sqrt{{\dot\gamma}^2_{2i-1}+{\dot\gamma}^2_{2i}}}{\vert s\kappa_{2i-1,2i}\vert},
\]
for $\kappa_{2i-1,2i}\ne 0$. Here $c_{2i-1}, c_{2i}, \varphi_{2i-1,2i}$ are the integration constants that depend on  the initial conditions. The period of rotation in the $\R^2(\gamma_{2i-1},\gamma_{2i})$-plane does not depend on the velocity:
\[
T_{2i-1,2i}=\frac{2\pi}{\vert \omega_{2i-1,2i}\vert}=\frac{2\pi m}{\vert s\kappa_{2i-1,2i}\vert}.
\]

In order to describe the first integrals of motion, we pass to the Lagrangian description. According to \eqref{AG}, we can take a Lagrangian with a term linear in velocities
\[
L_1^\Gamma(\gamma,\dot\gamma)=\frac{m}{2}\langle \dot\gamma,\dot\gamma\rangle+\frac{s}{2}\sum_{i=1}^{[n/2]} \kappa_{2i-1,2i}\big( (\gamma_{2i-1}+\Gamma_{2i-1})\dot\gamma_{2i}-(\gamma_{2i}+\Gamma_{2i})\dot\gamma_{2i-1}\big).
\]

The Lagrangian $L_1^\Gamma$ has the following Noether symmetries
\begin{equation}\label{simetrije}
\xi_i^{\Gamma_{2i-1},\Gamma_{2i}}=(\gamma_{2i-1}+\Gamma_{2i-1})\frac{\partial}{\partial \gamma_{2i}}-(\gamma_{2i}+\Gamma_{2i})\frac{\partial}{\partial \gamma_{2i-1}},
\end{equation}
that are also symmetries of the system with the Lagrangian $L$ modulo the gauge terms
\[
-\frac{d}{dt}s \frac{\kappa_{2i-1,2i}}{2}\big( (\gamma_{2i-1}+\Gamma_{2i-1})^2+(\gamma_{2i}+\Gamma_{2i})^2\big),\qquad i=1,\dots,[n/2].
\]

The corresponding gauge Noether first integrals of motion are  (e.g., see \cite{CS1981}):
\[
\langle \frac{\partial L}{\partial \dot\gamma},\xi_i^{\Gamma_{2i-1},\Gamma_{2i}}\rangle+s \frac{\kappa_{2i-1,2i}}{2}\big( (\gamma_{2i-1}+\Gamma_{2i-1})^2+(\gamma_{2i}+\Gamma_{2i})^2\big), \qquad i=1,\dots,[n/2].
\]
In the Hamiltonian description, they take the form
\[
\Phi_{2i-1,2i}^{\Gamma_{2i-1},\Gamma_{2i}}=(\gamma_{2i-1}+\Gamma_{2i-1})p_{2i}-(\gamma_{2i}+\Gamma_{2i}) p_{2i-1}+ s \frac{\kappa_{2i-1,2i}}{2}\big( (\gamma_{2i-1}+\Gamma_{2i-1})^2+(\gamma_{2i}+\Gamma_{2i})^2\big).
\]

Among the first integrals of motion $H_{2i-1,2i}, \Phi_{2i-1,2i}^{\Gamma_{2i-1},\Gamma_{2i}}$, there are three independent ones on
\[
(\R^4(\gamma_{2i-1},\gamma_{2i},p_{2i-1},p_{2i}),\Omega_{2i-1,2i}+\mathbf F_{2i-1,2i}).
\]

The systems \eqref{res1}, \eqref{res2}  are integrable in  non-commutative sense (see \cite{N, MF}). If $\kappa_{2i-1,2i}=0$, the trajectories
in $\R^4(\gamma_{2i-1},\gamma_{2i},p_{2i-1},p_{2i})$ are straight lines. Otherwise,  all the trajectories, outside $p_{2i-1}=p_{2i}=0$, are periodic with the period $T_{2i-1,2i}$. 
 As a result we have the following statement:

\begin{thm}\label{ocigledna}
For every $n\in \mathbb N$, the magnetic flows \eqref{Rn} are integrable in the  noncommutative sense.
Let $\kappa_{2i-1,2i}\ne 0$, $i=1,\dots,[n/2]$.
For even $n$, generic motions are quasi-periodic over $n/2$-dimensional isotropic tori $\mathbb T^{n/2}$, while for odd $n$, generic motions are quasi-periodic over isotropic cylinders $\mathbb T^{[n/2]}\times\R$.
\end{thm}

\begin{rem} For even $n$,
in the case when all $\vert\kappa_{2i-1,2i}\vert\ne 0$ are proportional:
\begin{equation*}
m_1\vert\kappa_{1,2}\vert^{-1}=\dots=m_{n/2}\vert\kappa_{n-1,n}\vert^{-1}, \qquad m_i\in\mathbb N,
\end{equation*}
such that the common divisor of $m_i$ is 1,
all the trajectories are closed with the same period
\[
T=\frac{m_1{2\pi m}}{{\vert s\kappa_{1,2}\vert}}=\dots=\frac{m_{n/2}{2\pi m}}{{\vert s\kappa_{n-1,n}\vert}}.
\]
The system is thus super-integrable.
\end{rem}

\section{Restrictions to  unit spheres}

In this section we consider  motion of a material point of mass $m$ and charge $-s$ restricted to the unit sphere $S^{n-1}\subset\R^n$ placed in a homogeneous magnetic field \eqref{polje}.
We work in redundant coordinates and consider the phase space $T^*S^{n-1}$ as a submanifold of $\mathbb{R}^{2n}(\gamma,p)$ given by the equations
\[
\phi_1=\langle\gamma,\gamma\rangle=1, \qquad \phi_2=\langle p,\gamma\rangle=0,
\]
and endowed with the twisted symplectic form $\omega+\mathbf f$, $\omega=\Omega\vert_{T^*S^{n-1}}$, $\mathbf f=\mathbf F\vert_{T^*S^{n-1}}$.
It is convenient to use the Dirac magnetic Poisson brackets on
\[
\R^{2n}_*=\{(\gamma,p)\in \mathbb{R}^{2n}\, \vert\, \phi_1\ne 0\}
\]
defined by (see e.g. \cite{AKN}):
\begin{equation}\label{eq:Dirac}
\{F, G\}_d =\{F, G\}^\kappa -\frac{\{F,\phi_1 \}^\kappa\{G,\phi_2\}^\kappa- \{F,\phi_2\}^\kappa
\{G,\phi_1\}^\kappa }{ \{\phi_1, \phi_2 \}^\kappa },
\end{equation}
where $\{\cdot, \cdot\}^\kappa$ are
the magnetic Poisson brackets on $\mathbb{R}^{2n}(\gamma,p)$ with respect to the twisted symplectic form $\Omega+\mathbf F$:
\[
\{F,G\}^\kappa=\sum_{i}\big(\frac{\partial
F}{\partial\gamma_i} \frac{\partial G}{\partial
p_i}-\frac{\partial F}{\partial p_i}\frac{\partial
G}{\partial\gamma_i}\big)+s\sum_{i,j} \kappa_{ij}\frac{\partial
F}{\partial p_i} \frac{\partial G}{\partial p_j}.
\]

The constraint functions  $\phi_1$ and $\phi_2$ are Casimir functions of the Dirac brackets \eqref{eq:Dirac}. The symplectic leaf
$\phi_1=1$, $\phi_2=0$ within $(\mathbb{R}^{2n}_*(\gamma,p), \{\cdot,\cdot\}_d)$ coincides with $(T^*S^{n-1},\omega+\mathbf f)$.

Let us consider the Hamiltonian flow on the Poisson manifold $(\mathbb{R}^{2n}_*(\gamma,p), \{\cdot,\cdot\}_d)$ defined by the Hamiltonian \eqref{HAM}:
\begin{align*}
\dot\gamma_i=& \{\gamma_i,H\}_d= \frac{1}{m} p_i-\frac{1}{m}\frac{\langle p,\gamma\rangle}{\langle \gamma,\gamma\rangle} \gamma_i,\\
\dot p_i   =&  \{p_i,H\}_d=\frac{s}{m}\langle \kappa p,\mathbf e_i\rangle+
\frac{{s}\langle p,\kappa\gamma \rangle-\langle p,p\rangle}{m\langle \gamma,\gamma\rangle }\gamma_i+\frac{1}{m}\frac{\langle p,\gamma\rangle}{\langle \gamma,\gamma\rangle} p_i -\frac{s}{m}\frac{\langle p,\gamma\rangle}{\langle \gamma,\gamma\rangle}\langle \kappa\gamma,\mathbf e_i\rangle.
\end{align*}

On the cotangent bundle $T^*S^{n-1}\subset \R^{2n}_*$ the equations of motion simplify to
\begin{equation}\label{magGF}
\dot\gamma= \frac{1}{m} p,\qquad
\dot p  =\frac{s}{m}\kappa p+\mu\gamma, \qquad \mu=\frac{1}{m}\big({s}\langle p,\kappa\gamma \rangle-\langle p,p\rangle\big).
\end{equation}
The function $\mu$ is the Lagrange multiplier and $\mu \gamma$ is the reaction force of the holonomic constraint $\phi_1=1$.

From now on, we consider a basis $[\mathbf e_1,\dots,\mathbf e_n]$ of $\R^n$ in which the magnetic form $\mathbf F$ takes the form \eqref{polje*}.
Thus, the system \eqref{magGF} becomes
\begin{align}
\label{eq1} &{\dot\gamma}_{2i-1}=\frac{1}{m}p_{2i-1},\qquad {\dot p}_{2i-1}=\frac{s}{m}\kappa_{2i-1,2i} p_{2i}+\mu\gamma_{2i-1},\\
\label{eq2} &{\dot\gamma}_{2i} =\frac{1}{m}p_{2i},\qquad\qquad  {\dot p}_{2i}=-\frac{s}{m}\kappa_{2i-1,2i} p_{2i-1}+\mu\gamma_{2i}, \qquad i=1,\dots,[n/2],
\end{align}
for $n$ even,
and, for $n$  odd, there is an additional equation:
\begin{equation}\label{eq3}
\dot\gamma_n=\frac{1}{m}p_n,\qquad\qquad {\dot p}_n=\mu\gamma_n.
\end{equation}

We can always perform a change of the basis with the permutations $\mathbf e_{2i-1} \longleftrightarrow \mathbf e_{2i}$ such that all parameters
$\kappa_{2i-1,2i}$ have the same sign, say
\begin{equation}\label{pozitivni}
\kappa_{2i-1,2i}\ge 0, \qquad i=1,\dots,[n/2].
\end{equation}
From now, on we assume that the inequalities \eqref{pozitivni} hold.

Note that the gauge Noether symmetries \eqref{simetrije} are tangent to the sphere $S^{n-1}$ for $\Gamma=0$,
leading to the following statement.

\begin{lem}\label{lema1}
The functions
\[
\Phi_{2i-1,2i}=\Phi_{2i-1,2i}^{0,0}=\gamma_{2i-1}p_{2i}-\gamma_{2i}p_{2i-1}+ s \frac{\kappa_{2i-1,2i}}{2}\big( \gamma_{2i-1}^2+\gamma_{2i}^2\big).
\]
are first integrals of motion of the magnetic flows \eqref{eq1} and \eqref{eq2} for even $n$ and   \eqref{eq1},  \eqref{eq2}, and   \eqref{eq3} for odd $n$. The first integrals of motion
Poisson commute: 
\[
\{\Phi_{2i-1,2i},\Phi_{2j-1,2j}\}_d=0, \qquad i,j=1,\dots,[n/2]
\]
on the Poisson manifold $(\mathbb{R}^{2n}_*(\gamma,p), \{\cdot,\cdot\}_d)$.
\end{lem}

We thus get the following theorem (see \cite{DGJ2023}):

\begin{thm}\label{stara}
Assume that at least one of the entries of the magnetic field \eqref{polje*}, e.g., $\kappa_{12}$ is not equal to zero. The magnetic  flows \eqref{magGF} are completely integrable on $S^2$ and $S^3$, corresponding to $n=3$ and $n=4$ respectively.
\end{thm}

\begin{rem}
Although the systems  \eqref{magGF} are quite natural and for $n=3$  well-known (see e.g.  \cite{S2002}), until recently a general form of the system \eqref{magGF} was not studied. In \cite{DGJ2023}, we derived equations \eqref{magGF} as the reduced equation of a gyroscopic nonholonomic Chaplygin system of a ball rolling over a sphere. The considered nonholonomic problem is a generalization of the Demchenko integrable system where the inertia operator of the system ball + gyroscope is proportional to the identity operator (see \cite{DGJ2020}). In \cite{DGJ2023}, we also
performed explicit integrations of the equations of motion \eqref{magGF} of the systems on $S^2$ and $S^3$, corresponding to $n=3$ and $n=4$ respectively, in elliptic functions. To aline the notation from this paper
with  \cite{DGJ2023}, one should set
\[
m={\tau}/{\varepsilon^2}, \qquad  s=1/{\varepsilon^2} \qquad   (s/m=1/\tau).
\]
\end{rem}

\begin{rem}
For even $n$, the one-form
\[
\mathbf a=\frac{s}{2}\sum_{i=1}^{n/2} \kappa_{2i-1,2i}(\gamma_{2i-1}d\gamma_{2i}-\gamma_{2i}d\gamma_{2i-1})\vert_{S^{n-1}}
\]
is a contact form on the sphere $S^{n-1}$ (see e.g. \cite{JJ2015}).
In the case when
\begin{equation}\label{jednakost}
\kappa_{12}=\kappa_{34}=\dots=\kappa_{n-1,n},
\end{equation}
the contact form $\mathbf a$ is proportional
to the \emph{standard contact form} on $S^{n-1}$. The magnetic geodesics on the standard contact sphere $S^{n-1}$ endowed with a metric deformed in the direction of the appropriate Reeb field were studied in \cite{DIMN2015}. For $n=4$, this system represents a magnetic flow on the Berger sphere \cite{IM2024}.
\end{rem}

In addition to the first integrals of motion $\Phi_{2i-1,2i}$, we construct  the following one.

\begin{lem}\label{lema2}
The function $J$ given by
\begin{align*}
J=& \frac{s^2}{m^2}\sum_{i=1}^{[n/2]}\kappa_{2i-1,2i}^2 (p_{2i-1}^2+p_{2i}^2)-\mu^2, \\
\mu=& \frac{s}{m}\sum_{i=1}^{[n/2]}\kappa_{2i-1,2i}(p_{2i-1}\gamma_{2i}-p_{2i}\gamma_{2i-1})-2H.
\end{align*}
is the first integral of motion of the magnetic flows  \eqref{eq1} and \eqref{eq2} for even $n$ and   \eqref{eq1},  \eqref{eq2}, and   \eqref{eq3} for odd $n$ on $T^*S^{n-1}$.
\end{lem}

\noindent\emph{Proof.}
We have $\dot H=0$ for $\phi_1=1$, $\phi_2=0$. Thus,
\begin{align*}
\dot J=& 2\frac{s^2}{m^2}\sum_{i=1}^{[n/2]}\kappa_{2i-1,2i}^2 \big(p_{2i-1}(\frac{s}{m}\kappa_{2i-1,2i} p_{2i}+\mu\gamma_{2i-1})
 +p_{2i}(-\frac{s}{m}\kappa_{2i-1,2i} p_{2i-1}+\mu\gamma_{2i})\big)\\
 &-2\mu\frac{s}{m}\sum_{i=1}^{[n/2]}\kappa_{2i-1,2i}\frac{d}{dt}\big(\gamma_{2i}p_{2i-1}-\gamma_{2i-1}p_{2i}\big)\\
 =&2\frac{s^2}{m^2}\sum_{i=1}^{[n/2]}\kappa_{2i-1,2i}^2 \big(\mu p_{2i-1}\gamma_{2i-1}+\mu p_{2i}\gamma_{2i})\big)+2\mu\frac{s}{m}\sum_{i=1}^{[n/2]}\kappa_{2i-1,2i}\frac{d}{dt}\big(\gamma_{2i-1}p_{2i}-\gamma_{2i}p_{2i-1}\big)
\end{align*}
On the other hand, by using the Noether first integrals $\Phi_{2i-1,2i}$, we get the identities
\[
\frac{d}{dt}\big(\gamma_{2i-1}p_{2i}-\gamma_{2i}p_{2i-1}\big)=- \frac{s\kappa_{{2i-1,2i}}}{m}\big( \gamma_{2i-1}p_{2i-1}+\gamma_{2i}p_{2i}\big), \qquad i=1,\dots,[n/2].
\]
Therefore, $\dot J=0$ for $\phi_1=1$, $\phi_2=0$.
\hfill $\Box$

\medskip

The proofs of the following two lemmas contain long calculations and we will omit them.

\begin{lem}\label{lema3}
The following commuting relations among the first integrals $J$, $ \Phi_{2i-1,2i}$, $ i=1,\dots,[n/2]$, take place:
\begin{align*}
\{J,\Phi_{2i-1,2i}\}_d=0, \qquad i=1,\dots,[n/2]
\end{align*}
on the Poisson manifold $(\mathbb{R}^{2n}_*(\gamma,p), \{\cdot,\cdot\}_d)$.
\end{lem}


\begin{lem}\label{lema4}
The functions $H$, $J$, $\Phi_{2i-1,2i}$, $i=1,\dots,[n/2]$ are functionally independent on $T^*S^{n-1}$ for $n\ge 5$ for all odd $n$ and all $\kappa$ and if  $n$ is even and $\kappa$ does not satisfy
$\kappa_{12}=\kappa_{34}=\dots=\kappa_{n-1,n}$.
\end{lem}

Let us  comment on the case when $n$ is even and the relation \eqref{jednakost} is satisfied.
Then the Lagrange multiplier $\mu$ is a first integral of motion,
\[
\mu= \frac{s\kappa_{12}}{m}\sum_{i=1}^{n/2}(p_{2i-1}\gamma_{2i}-p_{2i}\gamma_{2i-1})-2H=
\frac{s^2}{2m}\kappa^2_{12}-\frac{s}{m}\kappa_{12}\sum_{i=1}^{n/2}\Phi_{2i-1,2i}-2H,
\]
and
\begin{align*}
J=2\frac{s^2}{m}\kappa_{1,2}^2 H-\mu^2.
\end{align*}

It is  convenient to rewrite
the equations of motion \eqref{eq1} and \eqref{eq2} in a complex notation. Let us set
\[
z_i=\gamma_{2i-1}+\sqrt{-1}\gamma_i, \qquad w_i=p_{2i-1}+\sqrt{-1} p_{2i}, \qquad i=1,\dots,[n/2].
\]
Then the equations \eqref{eq1} and \eqref{eq2} take the form:
\begin{align}
\label{com1} \dot z_i=\frac{1}{m} w_i, \qquad \dot w_i=-\sqrt{-1}\frac{s\kappa_{2i-1,2i}}{m} w_i+\mu z_i, \qquad i=1,\dots,[n/2].
\end{align}

Thus,  the $SO(2)^{[n/2]}$--symmetry of the equations \eqref{eq1} and \eqref{eq2}
can be naturally considered as  the $U(1)^{[n/2]}$--symmetry of the equations \eqref{com1}.
If some of the parameters $\kappa_{2i-1,2i}$ are equal, then  there exists an additional $U(2)$--symmetry of the system, leading to the following first integrals of motion:

\begin{lem}\label{lema5}
If $\kappa_{2i-1,2i}=\kappa_{2j-1,2j}$, $i<j$, then the following functions
\begin{align*}
\Psi_{2i-1, 2i; 2j-1,2j}^1=&(\gamma_{2i}p_{2j-1}-\gamma_{2j-1}p_{2i})-(\gamma_{2i-1}p_{2j}-\gamma_{2j}p_{2i-1})\\
&-s\kappa_{2i-1,2i}(\gamma_{2i-1}\gamma_{2j-1}+\gamma_{2i}\gamma_{2j})\\
\Psi_{2i-1, 2i; 2j-1,2j}^2=&(\gamma_{2i-1}p_{2j-1}-\gamma_{2j-1}p_{2i-1})+(\gamma_{2i}p_{2j}-\gamma_{2j}p_{2i})\\
&-s\kappa_{2i-1,2i}(\gamma_{2i-1}\gamma_{2j}-\gamma_{2i}\gamma_{2j-1})
\end{align*}
are the first integrals of motion of the magnetic geodesic flow on $T^*S^{n-1}$. Moreover, the polynomials
\[
\Phi_{2i-1,2i},\, \Phi_{2j-1,2j},\, \Psi_{2i-1, 2i; 2j-1,2j}^1, \, \Psi_{2i-1, 2i; 2j-1,2j}^2
\]
generate the following four-dimensional Lie algebra
\begin{align*}
&\{\Phi_{2i-1,2i},\Phi_{2j-1,2j}\}_d=0,\\
&\{\Phi_{2i-1,2i}, \Psi_{2i-1, 2i; 2j-1,2j}^1\}_d=-\Psi_{2i-1, 2i; 2j-1,2j}^2,\\
&\{\Phi_{2j-1,2j}, \Psi_{2i-1, 2i; 2j-1,2j}^1\}_d=\Psi_{2i-1, 2i; 2j-1,2j}^2,\\
&\{\Phi_{2i-1,2i}, \Psi_{2i-1, 2i; 2j-1,2j}^2\}_d= \Psi_{2i-1, 2i; 2j-1,2j}^1,\\
&\{\Phi_{2j-1,2j}, \Psi_{2i-1, 2i; 2j-1,2j}^2\}_d=- \Psi_{2i-1, 2i; 2j-1,2j}^1,\\
&\{\Psi_{2i-1, 2i; 2j-1,2j}^1, \Psi_{2i-1, 2i; 2j-1,2j}^2\}_d=2\Phi_{2j-1,2j}-2\Phi_{2i-1,2i}.
\end{align*}
on the Poisson manifold $(\mathbb{R}^{2n}_*(\gamma,p), \{\cdot,\cdot\}_d)$. The Lie algebra is isomorphic to the reductive Lie algebra $\mathrm{so}(3)\oplus \mathbb{R}\cong \mathrm{u}(2)$.
\end{lem}

\noindent\emph{Proof.} The preservation of the polynomials $\Psi_{2i-1, 2i; 2j-1,2j}^1$, $\Psi_{2i-1, 2i; 2j-1,2j}^2$ along the magnetic flow \eqref{eq1}, \eqref{eq2} (and \eqref{eq3} for odd $n$) on $T^*S^{n-1}$ and a verification of the  Dirac-Poisson brackets follow by direct calculations.

Further, the function $e_0=\Phi_{2i-1,2i}+\Phi_{2j-1,2j}$ commutes with all the fisrt integrals. Let us set
\[
e_1=-\frac12\Psi_{2i-1, 2i; 2j-1,2j}^1,\quad e_2=-\frac12\Psi_{2i-1, 2i; 2j-1,2j}^2,\quad e_3=\frac{1}{2}(\Phi_{2j-1,2j}-\Phi_{2i-1,2i}).
\]
The Dirac-Poisson bracket   becomes
\[
\{e_1,e_2\}_d=e_3,\quad \{e_2,e_3\}_d=e_1,\quad \{e_3,e_1\}_d=e_2,\quad \{e_0,e_i\}_d=0,\quad i=1,2,3,
\]
which is the standard Lie bracket on $\mathrm{so}(3)\oplus \mathbb{R}\cong \mathrm{u}(2)$.
\hfill $\Box$

\medskip

Now we can state the main result of the paper.
Assume as before that at least one of the entries of the magnetic field \eqref{polje*}, e.g., $\kappa_{12}$ is not equal to zero.

\begin{thm}\label{glavna}
The magnetic geodesic flows are Liouville integrable on $T^*S^4$ and $T^*S^5$, corresponding to $n=5$ and $n=6$ respecively, for any $\kappa$. Moreover:

{\rm (i)} If $n=5$ and $\kappa_{12}=\kappa_{34}$, then the magnetic system \eqref{eq1}, \eqref{eq2},  and \eqref{eq3}  is integrable in the non-commutative sense: generic
motions are quasi-periodic over $3$-dimensional invariant
isotropic submanifolds.

{\rm (ii)} If $n=6$ and $\kappa_{12}=\kappa_{34}\ne\kappa_{56}$, then  the magnetic system \eqref{eq1} and  \eqref{eq2}  is integrable in the non-commutative sense:
generic motions are quasi-periodic over
$4$-dimensional invariant isotropic submanifolds.

{\rm (iii)} If $n=6$ and $\kappa_{12}=\kappa_{34}=\kappa_{56}$,  then  the magnetic system \eqref{eq1} and  \eqref{eq2}  is integrable in the non-commutative sense:
generic motions are quasi-periodic over
 $2$-dimensional invariant isotropic submanifolds.

{\rm (iv)}  For  $n=5$, $\kappa_{34}=0$ and for $n=6$, $\kappa_{34}=\kappa_{56}=0$, the magnetic  systems 
 are integrable in the non-commutative sense:
generic motions are quasi-periodic over
$3$-dimensional invariant isotropic submanifolds.
\end{thm}

\noindent\emph{Proof.}
If $\kappa_{12}\ne \kappa_{34}$ for $n=5$ or $\kappa_{12}\ne\kappa_{34}\ne\kappa_{56}$ for $n=6$, the systems are  Liouville integrabile, according to Lemmas \ref{lema1}, \ref{lema2}, \ref{lema3}, \ref{lema4}.

Let $n=5$ and $\kappa_{12}=\kappa_{34}$. The system is Liouville integrable as well. However,
there are following additional first integrals:
\begin{align*}
&\Psi_{12;34}^1=(\gamma_{2}p_3-\gamma_{3}p_2)-(\gamma_{1}p_4-\gamma_{4}p_1) -s\kappa_{12}(\gamma_1\gamma_3+\gamma_2\gamma_4),\\
&\Psi_{12;34}^2=(\gamma_{1}p_3-\gamma_{3}p_1)+(\gamma_{2}p_4-\gamma_{4}p_2) -s\kappa_{12}(\gamma_1\gamma_4-\gamma_2\gamma_3).
\end{align*}
Let $\mathcal F$ be the algebra of first integrals generated by $H,J,\Phi_{12},\Phi_{34}, \Psi_{12;34}^1, \Psi_{12;34}^2$.
By a  direct verification, we prove that the  first integrals $H,J,\Phi_{12},\Phi_{34}, \Psi_{12;34}^1$ are functionally independent. Thus, $\ddim \mathcal F=5$ (see \cite{BJ2003}). On the other hand, the functions
\[
H, \, J, \, \Phi_{12}+ \Phi_{34}
\]
Poisson commute with all the first integrals. Thus, $\dind\mathcal F=3$, $\ddim\mathcal F+\dind\mathcal F=\dim T^*S^4$,
and the system is completely integrable in the non-commutative sense with $3$-dimensional isotropic tori generated by the Hamiltonian
vector fields of $H$, $J$, and $\Phi_{12}+ \Phi_{34}$ (see \cite{N, BJ2003}).

Let $n=6$ and $\kappa_{12}=\kappa_{34}\ne \kappa_{56}$.
Again,  the Liouville integrability follows from Lemmas \ref{lema1}, \ref{lema2}, \ref{lema3}, \ref{lema4}.
According to Lemma \ref{lema5}, there are the following additional first integrals $\Psi_{12;34}^1$ and $\Psi_{12;34}^2$.
By a direct verification, we prove that the first integrals $H,J,\Phi_{12},\Phi_{34},\Phi_{56},\Psi_{12;34}^1$ are functionally independent.
The system is completely integrable in the non-commutative sense with $4$-dimensional isotropic tori, generated by the Hamiltonian
vector fields of $H$, $J$, $\Phi_{56}$, and $\Phi_{12}+ \Phi_{34}$.

Finally, let $n=6$ and $\kappa_{12}=\kappa_{34}=\kappa_{56}$.
Let $\mathcal F$ be the algebra of first integrals generated by
\[
H, \, \Phi_{12}, \, \Phi_{34}, \, \Phi_{56}, \, \Psi^1_{12;34}, \, \Psi^1_{12;56}, \, \Psi^1_{34;56}, \Psi^2_{12;34}, \, \Psi^2_{12;56}, \, \Psi^2_{34;56},
\]
 eight among them being independent. Since
$H$ and $\Phi_{12} + \Phi_{34} + \Phi_{56}$
commute with all first integrals, $\ddim\mathcal F+\dind\mathcal F=8+2=\ddim T^*S^5$ and the
system is completely integrable in the non-commutative sense with $2$-dimensional isotropic tori.
Note that the first integrals $\Phi_{12}, \, \Phi_{34}, \, \Phi_{56}, \, \Psi^1_{12;34}, \, \Psi^1_{12;56}, \, \Psi^1_{34;56}, \Psi^2_{12;34}, \, \Psi^2_{12;56}, \, \Psi^2_{34;56}$ form a Lie algebra isomorphic to $\mathrm{u}(3)$, which is related to the  $U(3)$-symmetry of the system \eqref{com1}.

The system is Liouville integrable as well:  $2$-dimensional isotropic tori can be non-uniquely organized
into $5$-dimensional Lagrangian tori, the level sets of the commuting first integrals (see \cite{BJ2003}). For example, one can take the following commuting first integrals
$H, \, \Phi_{12}, \, \Phi_{34}, \, \Phi_{56}, \, I$, where 
$I=2(\Phi_{12})^2+2(\Phi_{34})^2+(\Psi^1_{12;34})^2+(\Psi^2_{12;34})^2$.

The cases when $n=5$, $\kappa_{34}=0$ and $n=6$, $\kappa_{34}=\kappa_{56}=0$ are considered in Theorem \ref{integrabilni} below.
Note that item (iii) also follows from item (i) of Theorem \ref{integrabilni2}.
\hfill $\Box$

\begin{rem}
The first integral of motion $J+4H^2$ is a polynomial of the degree $3$ in momenta. Therefore, in Theorem \ref{glavna} the integrability is accomplished by means of  polynomial first integrals of the degrees $1$, $2$, and $3$ in momenta.
\end{rem}

\begin{rem}\label{n=4}
For $n=4$ and $\kappa_{12}=\kappa_{34}$, the elliptic integrals used in solving the equations become degenerate (see \cite{DGJ2023}).
The system is integrable in the non-commutative sense with generic invariant isotropic tori being  two-dimensional: there are four independent functions among  $H, \,\Phi_{12},\,\Phi_{34},\,\Psi^{1}_{12,34},\, \Psi^2_{12,34}$.
\end{rem}

\section{Reduction to the integrable cases with $n\le 6$}

For $n\ge 7$, by reduction to the cases with $n\le 6$, we derive the following three statements.

Similarly to Theorem 8.2 in \cite{DGJ2023}, we have:

\begin{thm}\label{integrabilni}
{\rm(i) $U(1)\times SO(n-2)$--symmetry.} Let $n\ge 5$. For $\kappa_{12}\ne 0$, $\kappa_{2i-1,2i}=0$, $i>1$ the magnetic geodesic flows  on  the sphere $S^{n-1}$ are completely integrable in noncommutative sense. Generic invariant manifolds are $3$-dimensional isotropic tori, the common level sets of $H$, $J$, $\Phi_{12}$, and the $SO(n-2)$--Noether first integrals
\[
\Phi_{ij}=\gamma_{i}p_{j}-\gamma_{j}p_{i}, \qquad 3\le i<j\le n.
\]
The magnetic flows are also Liouville integrable. One set of commuting first integrals consists of  $H$, $J$, $\Phi_{12}$, $\Phi_{34}$ and:
\begin{equation}\label{so(n-2)}
I_k=\sum_{3 \le i<j \le 4+k} (\Phi_{ij})^2, \qquad k=1,\dots, n-4.
\end{equation}

{\rm (ii) $U(2)\times SO(n-4)$--symmetry.} Let $n\ge 7$. For $\kappa_{12}=\kappa_{34} \ne 0$, $\kappa_{2i-1,2i}=0$, $i>2$ the magnetic geodesic flows  on  the sphere $S^{n-1}$ are completely integrable in the noncommutative sense. Generic invariant manifolds are $4$-dimensional isotropic tori, the common level sets of $H$, $J$,
$\Phi_{12}$, $\Phi_{34}$, $\Psi_{12;34}^1$, $\Psi_{12;34}^2$ and the $SO(n-4)$-Noether first integrals
\begin{align}\label{SO(n-4)}
\Phi_{ij}=\gamma_{i}p_{j}-\gamma_{j}p_{i}, \qquad 5\le i<j\le n.
\end{align}

{\rm (iii) $U(1)\times U(1)\times SO(n-4)$--symmetry.} Let $n\ge 7$. For $\kappa_{12}\ne \kappa_{34} \ne 0$, $\kappa_{2i-1,2i}=0$, $i>2$ the magnetic geodesic flow  on  the sphere $S^{n-1}$ is completely integrable in the noncommutative sense. Generic invariant manifolds are $5$-dimensional isotropic tori, the common level sets of $H$, $J$, $\Phi_{12}$, $\Phi_{34}$ and the $SO(n-4)$--Noether first integrals \eqref{SO(n-4)}.

In items {\rm (ii)} and {\rm (iii)} the magnetic flows are also Liouville integrable by means of commuting first integrals $H$, $J$, $\Phi_{12}$, $\Phi_{34}$, $\Phi_{56}$,
and:
\begin{align*}\label{so(n-4)}
I_k=\sum_{5 \le i<j \le 6+k} (\Phi_{ij})^2, \qquad k=1,\dots, n-6.
\end{align*}
\end{thm}

\noindent\emph{Sketch of the proof.} The essence of item (i) of Theorem \ref{integrabilni} is that for every trajectory $(\gamma(t),p(t))$ of the system, there exists a matrix $R\in SO(n-2)$ such that the trajectory
\[
\tilde\gamma(t)^T=\diag(1,1,R)\cdot \gamma(t)^T,\qquad \tilde p(t)^T=\diag(1,1,R)\cdot p(t)^T
\]
belongs to the invariant space
\[
\{(\gamma,p)\in\R^{2n}\, \vert\, \phi_1=1, \, \phi_2=0, \, \gamma_5=\dots=\gamma_n=0, \, p_5=\dots= p_n=0\}\subset T^* S^{n-1}
\]
diffeomorphic to $T^*S^3$ (see Theorem 8.2 in \cite{DGJ2023}). Therefore, the problem reduces to the problem of integrability of the magnetic system on $T^*S^3$ with $\kappa_{12}\ne 0$, $\kappa_{34}=0$ described in Theorem \ref{stara}.

The Poisson commuting integrals \eqref{so(n-2)} are related to the filtration of Lie algebras: $\mathrm{so}(2)<\mathrm{so}(3)<\dots<\mathrm{so}(n-2)$ (see e.g. \cite{JSV2023}).

The items (ii) and (iii) follow by a similar $SO(n-4)$--reduction to $T^*S^5$ and Theorem \ref{glavna}.
\hfill $\Box$

\medskip

Following the same line of reasoning we come to the following statement.

\begin{prop}\label{redukcija}
 Given $r\le [n/2]$, assume that
$
\kappa_{1,2}=\kappa_{3,4}=\dots=\kappa_{2r-1,2r}\ne 0.
$
Then:

{\rm (i)} The systems \eqref{com1} for $n$ even and the systems \eqref{com1}  and \eqref{eq3} for $n$  odd,  are $U(r)$--invariant, where the matrices $R$ from $U(r)$ act on  the first $r$ complex variables:
\[
(z_1,\dots,z_r)^T \longmapsto R\cdot (z_1,\dots,z_r)^T, \qquad  (w_1,\dots,w_r)^T \longmapsto R\cdot (w_1,\dots,w_r)^T
\]
 The corresponding first integrals $\Phi_{2i-1,2i}$, $\Psi_{2i-1,2i;2j-1,2j}^1$, $\Psi_{2i-1,2i;2j-1,2j}^2$, $1\le i \le r$, $i<j\le r$ form a Lie algebra
isomorphic to  $\mathrm{u}(r)$.

{\rm (ii)} In the real notation, for every trajectory $(\gamma(t),p(t))$ of the systems  \eqref{eq1} and  \eqref{eq2} for $n$ even and the systems
\eqref{eq1},  \eqref{eq2} and \eqref{eq3} for $n$ odd, there exists a matrix $R\in U(r)\subset SO(2r)$ such that the trajectory
\[
\tilde\gamma(t)^T=\diag(R,1,\dots,1)\cdot \gamma(t)^T, \qquad \tilde p(t)^T=\diag(R,1,\dots,1)\cdot p(t)^T
\]
belongs to the invariant space
\[
\{(\gamma,p)\in\R^{2n}\, \vert\, \phi_1=1, \, \phi_2=0, \, \gamma_1=\dots=\gamma_{2r-4}=0, \, p_1=\dots= p_{2r-4}=0\}\subset T^* S^{n-1}
\]
diffeomorphic to $T^*S^{n-2r+3}$.
\end{prop}

By combining  Theorem \ref{glavna}  and Proposition \ref{redukcija} we obtain

\begin{thm}\label{integrabilni2}
 Given $r\le [n/2]$, assume that
$\kappa_{12}=\kappa_{34}=\dots=\kappa_{2r-1,2r}\ne 0$. Then the following systems are integrable by means of  first integrals $H$, $J$,
$\Phi_{2k-1,2k}$, $k=1,\dots,[n/2]$, $\Psi_{2i-1,2i;2j-1,2j}^1$, $\Psi_{2i-1,2i;2j-1,2j}^2$, $1\le i \le r$, $i<j\le r$:

{\rm (i)}  When $n=2r$, the magnetic systems \eqref{eq1} and \eqref{eq2}  are integrable in the non-commutative sense: generic
motions are quasi-periodic over $2$-dimensional invariant isotropic submanifolds.

{\rm (ii)}  When $n=2r+1$,  the magnetic systems \eqref{eq1}, \eqref{eq2}, and \eqref{eq3}   are integrable in the non-commutative sense:
generic motions are quasi-periodic over $3$-dimensional invariant isotropic submanifolds.

{\rm (iii)}  When $n=2r+2$ and $\kappa_{12}\ne\kappa_{n-1,n}$,  the magnetic systems \eqref{eq1} and \eqref{eq2}   are integrable in the non-commutative sense: generic motions are quasi-periodic over $4$-dimensional invariant isotropic submanifolds.

The above systems are also Liouville integrable. One set of their commuting first integrals consists of
$H$, $J$, $\Phi_{2k-1,2k}$, $k=1,\dots,[n/2]$ and
\begin{equation}\label{komutativni}
I_k=\sum_{1 \le i < j \le k} (\Psi^1_{2i-1,2i;2j-1,2j})^2+(\Psi^2_{2i-1,2i;2j-1,2j})^2, \quad k=2,\dots,r.
\end{equation}
\end{thm}

Note that commuting first integrals are not unique (see \cite{BJ2003}) and the first integrals \eqref{komutativni} are related to the  filtration of Lie algebras $\mathrm{u}(1)<\mathrm{u}(2)<\dots<\mathrm{u}(n)$
(see e.g. \cite{JSV2023}).

\begin{rem}
Item (i) of Theorem \ref{integrabilni2} and item (ii) of Proposition \ref{redukcija} are in agreement with the
result of \cite{DIMN2015} for the standard contact sphere where  $n$ is even and $\kappa_{12}=\kappa_{34}=\dots=\kappa_{n-1,n}$.
There, it is proved that for every magnetic geodesic line $\gamma(t)$ on the standard contact sphere $S^{2r-1}$, there exists
a totally geodesic 3-sphere $S^3_{0}\subset S^{2r-1}$, such that $\gamma(t)$ is the magnetic geodesic line on $S^3_0$ endowed with the
standard contact structure.
\end{rem}

Finally,  there is a
 $U(r)\times SO(n-2r)$--symmetry
that allows a  reduction to $T^*S^5$. By combining
Theorems \ref{integrabilni} and \ref{integrabilni2}, we get:

\begin{thm}\label{integrabilni3}
 Given $r<[n/2]$, assume that
\[
\kappa_{12}=\dots=\kappa_{2r-1,2r}\ne 0,  \qquad \kappa_{2r+1,2r+2}=\dots=\kappa_{2[n/2]-1,2[n/2]}=0.
\]
Then the magnetic systems \eqref{eq1}, \eqref{eq2} for $n$ even and the magnetic  systems   \eqref{eq1}, \eqref{eq2} and \eqref{eq3} for $n$  odd, are integrable in the non-commutative sense by means of  first integrals $H$, $J$,
$\Phi_{2k-1,2k}$, $k=1,\dots,r$, $\Psi_{2i-1,2i;2j-1,2j}^1$, $\Psi_{2i-1,2i;2j-1,2j}^2$, $1\le i \le r$, $i<j\le r$ and
\[
\Phi_{ij}=\gamma_{i}p_{j}-\gamma_{j}p_{i}, \qquad 2r+1\le i<j\le n.
\]
Generic motions are quasi-periodic over $4$-dimensional invariant isotropic submanifolds.
The systems  are also Liouville integrable. One set of their commuting first integrals consists of
$H$, $J$, $\Phi_{2k-1,2k}$, $k=1,\dots,r$ and
\begin{align*}
& I_k=\sum_{1 \le i < j \le k} (\Psi^1_{2i-1,2i;2j-1,2j})^2+(\Psi^2_{2i-1,2i;2j-1,2j})^2, \quad k=2,\dots,r,\\
& I_{r+k}=\sum_{2r+1 \le i<j \le 2r+2+k} (\Phi_{ij})^2, \qquad k=1,\dots, n-2r-2.
\end{align*}
\end{thm}

The following classical mechanical problem is well-known:  motion of a material point of mass $m$ in $\R^n$ under the influence of the potential force field with a quadratic potential
\[
V(\gamma)=\frac12\big(a_1 \gamma_1^2+\dots+a_n\gamma_n^2\big).
\]
While the system of uncoupled harmonic oscillators is trivially integrable, the restriction of the problem to the sphere $S^{n-1}$ leads to one of the most interesting finite-dimensional  integrable systems  - the \emph{Neumann system} (see e.g. \cite{M}). The above considerations suggest  a similar situation in the case of the considered magnetic flows. There is an obvious integrability in $\R^n$, as shown in Theorem \ref{ocigledna}, while the integrability  of the system  restricted to the sphere $S^{n-1}$ is highly non-trivial.

\begin{con}\label{hipoteza}
Motion of a material point of mass $m$  on the  sphere $S^{n-1}$ placed in a constant homogeneous magnetic field \eqref{polje} is completely integrable in the Liouville or in the noncommutative sense for any $n$ and $\kappa$.
\end{con}

\begin{rem}
It is of interest to consider magnetic geodesic flows for  more general metrics on the spheres $S^{n-1}$.
For example, there is a well-known Volterra--Zhukovskii integrable case of motion of a rigid body with a gyroscope around a fixed point \cite{Vol, Zh}.
By lifting the rigid body metric with respect to the two-covering $S^3\to SO(3)$, we obtain an integrable magnetic geodesic flow for any left-invariant metric on $S^3$.
Next, one can apply left $SO(2)$-magnetic reduction (see \cite{KNP2005}) to get a suitable integrable magnetic flow on the sphere $S^2$.
Let us note that there is a difference between the nonholonomic gyroscopic reduction \cite{DGJ2023, GMD2024} and the Hamiltonian magnetic reduction \cite{KNP2005}.
\end{rem}

\subsection*{Acknowledgements} This work is dedicated to Academician V. V. Kozlov on the occasion of his 75-th anniversary. The authors have had pleasure to learn a lot from numerous  personal scientific interactions with Valery Vasil'evich and his students  as well as from their books, papers, and lectures.

We are grateful to the referees for useful remarks. This research was supported by the Serbian Ministry of Science, Technological Development and Innovation through Mathematical Institute of Serbian Academy of Sciences and Arts and the Simons Foundation grant no. 854861.


\begin{thebibliography}{99}

\small

\bibitem{AKN} V. I. Arnold, V. V. Kozlov, A. I. Neishtadt,
\emph{Mathematical Aspects of Classical and Celestial Mechanics},
Encylopaedia of Mathematical Sciences \textbf{3}, Springer-Verlag, Berlin, 1989.

\bibitem{AS2020} A. Arvanitoyeorgos, N. P. Souris, \emph{Motion of charged particle in a class of homogeneous spaces}, Math.
Phys. Anal. Geom. \textbf{23} (2020) 22, 8 pp.

\bibitem{BK2017}
S. V. Bolotin, V.  V. Kozlov, \emph{Topology, singularities and integrability in Hamiltonian systems with two degrees of freedom}, Izv. RAN. Ser. Mat., \textbf{81}(2017), no. 4, 3--19 (Russian); English translation: Izv. Math., \textbf{81}(2017) 671--687.

\bibitem{BJ2003} A. V. Bolsinov,  B. Jovanovi\' c,
\emph{Non-commutative integrability, moment map and geodesic flows}.
Annals of Global Analysis and Geometry {\bf 23} (2003) 305--322, arXiv: math-ph/0109031.

\bibitem{BJ2006} A. V. Bolsinov, B. Jovanovi\' c, \emph{Magnetic Geodesic Flows on Coadjoint Orbits}, J. Phys. A: Math. Gen.  {\bf 39} (2006) L247--L252,
arXiv: math-ph/0602016.

\bibitem{BJ2008} A. V. Bolsinov, B. Jovanovi\' c, \emph{Magnetic Flows on Homogeneous Spaces}, Comm. Math. Helv, \textbf{83} (2008) 679--700.


\bibitem{CS1981}
F. Cantrjin, W. Sarlet, \emph{Generalizations of Noether’s theorem in classical mechanics}, SIAM
Rev. \textbf{23} (1981), 467--494.


\bibitem{DGJ2020} V. Dragovi\' c, B. Gaji\' c, B. Jovanovi\' c,
\emph{Demchenko's nonholonomic case of a gyroscopic ball rolling without sliding over a sphere after his 1923 Belgrade doctoral thesis},
Theor. Appl. Mech. \textbf{47}(2) (2020), 257--287,   	arXiv:2011.03866.

\bibitem{DGJ2023} V. Dragovi\' c, B. Gaji\' c, B. Jovanovi\' c, \emph{Gyroscopic Chaplygin systems and integrable magnetic flows on spheres},
J. Nonlinear Sci. \textbf{33}(2023) 43,  arXiv:2110.09938 [math-ph].

\bibitem{DIMN2015} S. L. Druta-Romaniuc, J. Inoguchi, M. I. Munteanu,  A. I. Nistor, \emph{Magnetic curves
in Sasakian manifolds}, J. Nonlinear Math. Physics, \textbf{22} (2015) 3, 428--447.

\bibitem{E2005} D. I. Efimov,  \emph{The magnetic geodesic flow on a homogeneous symplectic manifold},
Sib. Mat. Zh. \textbf{46} (2005) 106–-118 (Russian). English translation:
Siberian Mathematical Journal {\bf 46} (2005) 83--93.

\bibitem{F} R. Feynman, \emph{The Motion of Charges in Electric and Magnetic Fields}, {Lectures in Physics}, \url{https://www.feynmanlectures.caltech.edu/II_29.html}.

\bibitem{GMD2024}
L. C. Garcia-Naranjo, J. C. Marrero, D. Martin de Diego, P. E. P. Valdes, \emph{Almost-Poisson Brackets for Nonholonomic Systems with Gyroscopic Terms and Hamiltonisation}, J. Nonlinear Sci. \textbf{34} (2024) 110 (52pp),  	arXiv:2309.11597 [math-ph].


\bibitem{IM2024} J. Inoguchi, M. I. Munteanu,  \emph{Homogeneity of magnetic trajectories in the Berger sphere}, arXiv:2406.15886 [math.DG].

\bibitem{JJ2015} B. Jovanovi\'c, V. Jovanovi\'c, \emph{Contact flows and integrable systems}, \textbf{87} (2015) 217--232, arXiv:1212.2918 [math.SG].


\bibitem{JSV2023} B. Jovanovi\' c, T. \v Sukilovi\' c, S. Vukmirovi\' c,  \emph{Integrable systems associated to the filtrations of Lie algebras}, Regul. Chaot. Dyn. \textbf{28} (2023), 44--61, arXiv:1912.03199 [nlin.SI].

\bibitem{Ko1997} V. V. Kozlov, \emph{Stabilization of the unstable equilibria of charges by intense magnetic fields}, Journal of Applied Mathematics and Mechanics
\textbf{61} (1997) Issue 3,  377--384.


\bibitem{KP} V. V. Kozlov, S. A. Polikarpov, \emph{
Periodic billiard trajectories in a magnetic field}, Journal of Applied Mathematics and Mechanics
\textbf{69} (2005), Issue 6,  844--851.


\bibitem{KNP2005}
 N. Kowalzig, N. Neumaier, M. J. Pflaum,
\emph{Phase Space Reduction of Star Products on Cotangent Bundles},
Annales Henri Poincaré \textbf{6} (2005) 485--552, arXiv:math/0403239 [math.SG]

\bibitem{MSY2008}
A. \,A. Magazev, I. \,V. Shirokov, Yu. \,A. Yurevich, \emph{Integrable magnetic geodesic flows on Lie groups},
TMF, \textbf{156} (2008), no 2,  189--206 (in Russian); English. transl. Theoret. and Math. Phys.,  \textbf{156} (2008) 1127--1141,
arXiv:1111.0726 [math-ph].


\bibitem{MF}
A. S. Mishchenko, A. T. Fomenko, \emph{Generalized Liouville
method of integration of Hamiltonian systems}. Funkts. Anal. Prilozh. {\bf 12} (1978) No.2, 46-56   (Russian); English translation:  Funct. Anal. Appl. {\bf
12}, (1978) No. 2,  113--121

\bibitem{M}
J. Moser,  \emph{Various aspects of integrable Hamiltonian systems}, Dynamical systems (C.I.M.E. Summer
School, Bressanone, 1978). In: Progr. Math. vol. 8, pp. 233–289. Birkhuser, Boston (1980)

\bibitem{N}
N. N. Nekhoroshev, \emph{Action-angle variables and their generalization}, Tr. Mosk. Mat. O.-va. \textbf{26}
(1972), 181--198 (Russian). English translation: {Trans. Mosc. Math. Soc.} {\bf 26} (1972)  180--198.

\bibitem{NS1981} S. P. Novikov, I. Shmel'tser, \emph{Periodic solutions
to Kirchoff equations for a free motion of a rigid body in fluid
and Lusternik-Shnirelman-Morse extended theory}. Funkc. Anal. Pril. {\bf 15} (1981) No. 3,  54--66 (Russian).
English translation: Funct. Anal. Appl. \textbf{15}(1981), No. 3, 197–-207

\bibitem{N1982}
S.P. Novikov,
\emph{The Hamiltonian formalism and a many-valued analogue of Morse theory}, UMN \textbf{37} (1982), No. 5, 3-–49 (Russian);
English translation: Russian Mathematical Surveys, \textbf{37} (1982), No. 5, 1–-56.


\bibitem{S2002}
P. Saksida,  \emph{Neumann system, spherical pendulum and magnetic fields}, J. Phys. A: Math. Gen. \textbf{35} (2002) 5237--5253.

\bibitem{T2016} I. \,A. Taimanov, \emph{On first integrals of geodesic flows on a two-torus},
Modern problems of mechanics, Collected papers, Tr. Mat. Inst. Steklova, \textbf{295}, MAIK Nauka/Interperiodica, Moscow, 2016, 241--260 (Russian);
English transl.: Proc. Steklov Inst. Math., \textbf{295} (2016) 225–242.


\bibitem{VoKo}
I. V. Volovich,  V. V. Kozlov, \emph{On Maxwell’s Equations with a Magnetic Monopole on Manifolds},  Proceedings of the Steklov Institute of Mathematics,
Volume 306,  (2019), 43--46.

\bibitem{Vol} V. Volterra, \emph{Sur la th\'eorie des variations des latitudes}, Acta Math. \textbf{22} (1899), no 1, 201--357 (French).



\bibitem{Zh} N. Ye. Zhukovskii, \emph{On the motion of a rigid body with cavities filled with a homogeneous liquid drop}, Collected Works Vol. 1, Gostekhizdat, Moscow, 1948, pp. 31--152 (Russian).

\end{thebibliography}
\end{document}